\documentclass[a4paper,12pt]{article}
\usepackage {amssymb,latexsym}
\usepackage {amsmath}
\usepackage[T1]{fontenc}
\usepackage{verbatim}
\usepackage{relsize}
\usepackage{graphicx}

\def\R{\mathbb{R}}

  \def\Aut{{\rm Aut}}

\usepackage{xypic}
\xyoption{all}
\input{xypic}
\newdir{ >}{{}*!/-8pt/\dir{>}}

\newtheorem{theorem}{Theorem}[section]
\newtheorem{prop}{Proposition}[section]

\newtheorem{cor}{Corollary}[section]
\newtheorem{defi}{Definition}[section]

                                                     \medskip
\newenvironment{prooof}{
                        \noindent{\bf\small Proof: }\small}
                                       {\hfill {$\mathbf \Box$}\\}

\makeatletter
\@addtoreset{equation}{section}

\makeatother

\begin{document}

\thispagestyle{empty}
\pagestyle{empty}
\title
   {\LARGE \bf A dirty integration of
       Leibniz algebras}


\author{ Martin Bordemann \\
Universit\'{e} de Haute Alsace, Mulhouse\\
  \texttt{Martin.Bordemann@uha.fr}\\
\and  Friedrich Wagemann\\
     Universit\'{e} de Nantes\\
     \texttt{Friedrich.Wagemann@univ-nantes.fr}}


\maketitle
\thispagestyle{empty}

\vspace{1cm}
\footnotesize
\noindent {\bf Abstract}\\[5mm]
\thispagestyle{empty}
\begin{minipage}{13cm}
 In this paper we present an integration of any real
 finite-dimensional Leibniz algebra as a Lie rack which reduces
 in the particular case of a Lie algebra
 to the ordinary connected simply connected Lie group.
 The construction is not functorial.
\end{minipage}
\normalsize

\newpage

\tableofcontents
 \thispagestyle{empty}
\newpage

\pagestyle{plain}

\section{Introduction}

All manifolds considered in this manuscript are assumed
to be Hausdorff and second countable.

Recall that a \emph{pointed rack} is a pointed set $(X,e)$ together with a 
binary operation 
$\rhd:X\times X\to X$ such that for all $x\in X$, the map $y\mapsto x\rhd y$ 
is bijective and such that for all $x,y,z\in X$,
the self-distributivity and unit relations
$$
   x\rhd(y\rhd z)\,=\,(x\rhd y)\rhd(x\rhd z),~~~e\triangleright x = x,~~~
   \mathrm{and}~~~x\triangleright e = e
$$
are satisfied. 
Imitating the notion of a Lie group, the smooth version of a pointed rack 
is called
{\it Lie rack}.  

An important class of examples of racks are the so-called 
{\it augmented racks}, 
see \cite{FR92}. 
An augmented rack is the data of a group $G$, a $G$-set $X$ and a map 
$p:X\to G$ 
such that for all $x\in X$ and all $g\in G$, 
$$p(g\cdot x)\,=\,gp(x)g^{-1}.$$
The set $X$ becomes then a rack by setting $x\rhd y\,:=\,p(x)\cdot y$. 
  
Lie racks are intimately related to
{\it Leibniz algebras} ${\mathfrak h}$, i.e. a vector space 
${\mathfrak h}$ with a 
bilinear bracket $[,]:{\mathfrak h}\otimes{\mathfrak h}\to{\mathfrak h}$ 
such that
for all $X,Y,Z\in{\mathfrak h}$, $[X,-]$ acts as a derivation:
$$[X,[Y,Z]]\,=\,[[X,Y],Z]+[Y,[X,Z]].$$  
Indeed, Kinyon showed in \cite{Kin07} that the tangent space at 
$e\in H$ of a Lie rack $H$ carries 
a natural structure of a Leibniz algebra, generalizing the relation 
between a Lie group 
and its tangent Lie algebra. Conversely, every 
(finite dimensional real or complex) 
Leibniz algebra ${\mathfrak h}$
may be integrated into a Lie rack (with underlying manifold ${\mathfrak h}$) 
using the rack product
\begin{equation}   \label{rack_structure}
X\rhd Y\,:=\,e^{{\rm ad}_X}(Y),
\end{equation}
noting that the exponential of the inner derivation ${\rm ad}_X$ for 
each $X\in{\mathfrak h}$
is an automorphism. Although the assignment 
$\big(\mathfrak{h},[~,~]\big)\to \big(\mathfrak{h},0, \triangleright \big)$
is functorial since morphisms of Leibniz algebras are easily seen to
go to morphisms of pointed Lie racks, the restriction to the category of 
all Lie 
algebras would not give the usual integration as a Lie group.

The purpose of the present paper is to construct an integration procedure 
which integrates real, finite-dimension Leibniz algebras into Lie racks in 
such a way that the restriction to Lie algebras gives the conjugation rack 
underlying the simply connected
Lie group corresponding to a (real, finite-dimensional) Lie algebra.

This problem has been encountered by J.-L. Loday in 1993 \cite{Lod} 
in the search of quantifying
the lack of periodicity in algebraic K-Theory. Several attempts and 
constructions have been published since then. In 2010, Simon Covez \cite{Cov}
solves 
in his thesis the local integration problem by constructing a local Lie 
rack integrating a given (real, finite dimensional) Leibniz algebra in such 
a way that in the case of Lie algebras, one obtains 
the conjugation Lie rack 
underlying the usual (simply connected) Lie group integrating it. Other 
important contributions to the problem include Mostovoy's article \cite{Mos} 
where he solves the problem in the framework of formal groups. 
The general problems is still open to our knowledge and our article 
constitutes another step towards its solution. 

In this article, we construct a (global) Lie rack integrating a given (real, 
finite dimensional) Leibniz algebra ${\mathfrak h}$ in such a way that in the 
case of a Lie 
algebra, the construction yields the conjugation rack underlying the 
usual (simply connected) Lie group 
integrating it. More precisely, we work with augmented Leibniz algebras, i.e. 
Leibniz algebras ${\mathfrak h}$ with an action of a 
Lie algebra ${\mathfrak g}$ by derivations and an equivariant map
$p:{\mathfrak h}\to{\mathfrak g}$ to ${\mathfrak g}$. 
We integrate the quotient Lie algebra $p({\mathfrak h})=:{\mathfrak g}'$ 
into a Lie group $G'$ and integrate its action 
on ${\mathfrak h}$ such that the resulting (global) augmented Lie rack is 
an affine bundle over $G'$ with typical fiber ${\mathfrak z}:=\mathrm{Ker}(p)$.
This is the content of our main theorem, Theorem \ref{main_theorem}. 

Perhaps the most interesting point of the article is the fact that for this 
construction, we need an open neighbourhood in the Lie group $G'$ on which 
the exponential is a diffeomorphism and which is invariant under (the 
connected component of the identity) $\Aut_0(G')$. It is not elementary
to show that such a neighbourhood exists, and the proof of it will take 
the 10 pages of Appendix A. The arguments are rather classical and are inspired 
by Lazard-Tits \cite{LT66} and \DJ okovi\'{c}-Hofmann \cite{DH97}.
On the other hand, it is the use of this neighbourhood which renders our 
construction non-functorial, which is a major drawback of the theory. 
This is why we call our integration {\it dirty}. For the moment, 
we do not know whether there exists a functorial construction of a Lie 
rack integrating a 
given Leibniz algebra (such that in the special case of a Lie algebra, we get 
back the conjugation rack underlying the usual simply connected Lie group).

\subsubsection*{Acknowledgements}
FW thanks Universit\'e de Haute Alsace (Mulhouse, France)  
for several research visits where 
the subject of this article has been discussed.

\section{Augmented Leibniz algebras and Lie racks}

\subsection{(Augmented) Leibniz algebras}
Let $K$ be any unital commutative ring containing the rational numbers.
We are mainly interested in the case $K=\R$.    
All modules in this section will be considered over $K$.

Recall that a \emph{Leibniz algebra} over $K$ is a $K$-module
$\mathfrak{h}$ equipped with a linear map $[~,~]:
\mathfrak{h}\otimes \mathfrak{h}\to\mathfrak{h}$, written
$x\otimes y\mapsto [x,y]$ such that the \emph{left Leibniz identity} holds
for all $x,y,z\in \mathfrak{h}$
 \begin{equation}\label{EqDefLeibnizIdentity}
  \big[x, [y, z]\big]=
     \big[ [x, y], z\big]
         + \big[y,[x,z]\big]
 \end{equation}
 A morphism of Leibniz algebras $\Phi:\mathfrak{h}\to \mathfrak{h}'$
 is a $K$-linear map preserving brackets, i.e. 
 for all $x,y\in\mathfrak{h}$ we have 
 $\Phi\big([x,y]\big)=\big[\Phi(x),\Phi(y)\big]'$.
Recall first that each Lie algebra over $K$ is a Leibniz algebra
giving rise to a functor $\mathsf{i}$ 
from the category of all Lie algebras (over $K$), 
$K\mathbf{LieAlg}$, to the 
category of all Leibniz algebras (over $K$),  $K\mathbf{Leib}$.

Furthermore, recall that each Leibniz algebra has two canonical 
$K$-submod\-ules
\begin{eqnarray}
   Q(\mathfrak{h})&:=&
   \big\{x\in\mathfrak{h}~|~\exists~N\in\mathbb{N}\setminus\{0\},~
    \exists~\lambda_1,\ldots,\lambda_N\in K,~\exists~x_1,\ldots,x_N~
    \nonumber\\
    & &
   ~~~~~~~~~~~~~~~~~~ 
       \mathrm{such~that~}x=\sum_{r=1}^N\lambda_r[x_r,x_r]\big\}, \\
    \mathfrak{z}(\mathfrak{h})
    & := & \big\{x\in\mathfrak{h}~|~\forall~y\in\mathfrak{h}
                     :~[x,y]=0\big\}.
\end{eqnarray}
It is well-known and not hard to deduce from the Leibniz identity 
that both
$Q(\mathfrak{h})$ and $\mathfrak{z}(\mathfrak{h})$ are two-sided
abelian ideals of $(\mathfrak{h}, [~,~])$, that
$Q(\mathfrak{h})\subset \mathfrak{z}(\mathfrak{h})$, and that
the quotient Leibniz algebras
\begin{equation}
    \overline{\mathfrak{h}}:=\mathfrak{h}/Q(\mathfrak{h})
    ~~~\mathrm{and}~~~
       \mathfrak{h}/\mathfrak{z}(\mathfrak{h})
\end{equation}
are Lie algebras. Since the ideal $Q(\mathfrak{h})$ is clearly
mapped
into the ideal $Q(\mathfrak{h}')$ by any morphism of Leibniz algebras
$\mathfrak{h}\to\mathfrak{h}'$ (which is a priori not the case
for $\mathfrak{z}(\mathfrak{h})$ !), there is an obvious functor
$\mathfrak{h}\to \overline{\mathfrak{h}}$ from the category of
all Leibniz algebras to the category of all Lie algebras.
It is not hard to see and not important for the sequel that the functor
$\mathfrak{h}\to \overline{\mathfrak{h}}$ is a left adjoint functor of
the inclusion functor of the category of all Lie algebras in the category of
all Leibniz algebras whence the former is a reflective subcategory of the latter,
see e.g.~\cite[p.91]{Mac98} for definitions.

It is easy to observe that in both cases of the above Lie algebras,
$\overline{\mathfrak{h}}$ and $\mathfrak{h}/\mathfrak{z}(\mathfrak{h})$, 
there is the following structure:
\begin{defi}\label{DAugementedLeibnizAlgebra}
  A quintuple 
  $\big(\mathfrak{h},p,\mathfrak{g},[~,~]_\mathfrak{g},\dot{\rho}\big)$
  is called a  $\mathfrak{g}$-\emph{augmented Leibniz algebra} iff
  the following holds:
  \begin{enumerate}
   \item $\big(\mathfrak{g},[~,~]_\mathfrak{g}\big)$ is a Lie algebra over $K$.
   \item $\mathfrak{h}$ is a $K$-module which is a left $\mathfrak{g}$-module
     via the $K$-linear map $\dot{\rho}:\mathfrak{g}\otimes \mathfrak{h}\to
     \mathfrak{h}$ written $\dot{\rho}_\xi(x)=\xi.x$ for all $\xi\in\mathfrak{g}$
     and $x\in\mathfrak{h}$.
    \item $p:\mathfrak{h}\to \mathfrak{g}$ is a $K$-linear morphism of 
        $\mathfrak{g}$-modules, i.e.~for all $\xi\in\mathfrak{g}$
     and $x\in\mathfrak{h}$
       \begin{equation}\label{EqDefPhiMorphismOfGLieModules}
          p(\xi.x)=\big[\xi, p(x)\big]_\mathfrak{g}.
       \end{equation}
   \end{enumerate}
   A \emph{morphism of augmented Leibniz algebras} 
   $\big(\mathfrak{h},p,\mathfrak{g},[~,~]_\mathfrak{g},\dot{\rho}\big)$
   $\to$ $\big(\mathfrak{h}',p',\mathfrak{g}',[~,~]'_\mathfrak{g},
   \dot{\rho}'\big)$ is a pair $(\Phi,\phi)$ of $K$-linear maps where
   $\phi:\mathfrak{g}\to\mathfrak{g}'$ is a morphism of Lie algebras,
   $\Phi:\mathfrak{h}\to\mathfrak{h}$ is a morphism of Lie algebra modules over
   $\phi$, i.e. ~for all $x\in\mathfrak{h}$ and $\xi\in\mathfrak{g}$
   \begin{equation}\label{EqDefPhiphiAugLeib}
      \Phi(\xi.x)=\phi(\xi).\Phi(x).
   \end{equation}
   Moreover the obvious diagram commutes, i.e.
   \begin{equation}\label{EqDefPhiphiPPPrimeAugLeib}
      p'\circ \Phi = \phi \circ p.
   \end{equation}
\end{defi}

The following properties are immediate from the definitions:
\begin{prop}
  Let $\big(\mathfrak{h},p,\mathfrak{g},[~,~]_\mathfrak{g},\dot{\rho}\big)$
  be an augmented Leibniz algebra. Define the following bracket on 
  $\mathfrak{h}$:
    \begin{equation}\label{EqDefLeibnizBracketForAugmentedLeibnizAlg}
          [x,y]_\mathfrak{h} := p(x).y.
    \end{equation}
  \begin{enumerate}
   \item $\big(\mathfrak{h},[,]_\mathfrak{h}\big)$ is a Leibniz algebra
      on which $\mathfrak{g}$ acts as derivations. If $(\Phi,\phi)$
      is a morphism of augmented Leibniz algebras, then
      $\Phi$ is a morphism of Leibniz algebras.
   \item The kernel of $p$, $\mathrm{Ker}(p)$, is a $\mathfrak{g}$-invariant
     two-sided abelian ideal of $\mathfrak{h}$ satisfying
     $Q(\mathfrak{h})\subset \mathrm{Ker}(p) \subset \mathfrak{z}(\mathfrak{h})$.
   \item The image of $p$, $\mathrm{Im}(p)$, is an ideal of the Lie algebra
      $\mathfrak{g}$.
  \end{enumerate}
\end{prop}
\begin{prooof}
 We just check the Leibniz identity: Let $x,y,z\in \mathfrak{h}$, then,
 writing $[~,~]_\mathfrak{h}=[~,~]$,
 \begin{eqnarray*}
    \big[x, [y, z]\big] & = &  p(x).\big(p(y).z\big) \\
      &=&  p(x).\big(p(y).z\big) - p(y).\big(p(x).z\big) + 
              p(y).\big(p(x).z\big)\\
      & = & \big[p(x),p(y)\big]_\mathfrak{g}.z + \big[y, [x, z]\big] 
       \stackrel{(\ref{EqDefPhiMorphismOfGLieModules})}{=}  
       p\big(p(x).y\big).z + \big[y, [x, z]\big] \\
       & = & \big[[x,y],z\big]  + \big[y, [x, z]\big].
 \end{eqnarray*}
\end{prooof}

It follows that the class of all augmented Leibniz algebras forms a category
$K\mathbf{LeibA}$, and there is an obvious forgetful functor from 
$K\mathbf{LeibA}$
to $K\mathbf{Leib}$ associating to 
$\big(\mathfrak{h},p,\mathfrak{g},[~,~]_\mathfrak{g},\dot{\rho}\big)$
the Leibniz algebra $\big(\mathfrak{h}, [~,~]_\mathfrak{h}\big)$ where the 
Leibniz bracket $[~,~]_\mathfrak{h}$ is defined in eqn 
(\ref{EqDefLeibnizBracketForAugmentedLeibnizAlg}).

On the other hand there is a functor from $K\mathbf{Leib}$ to
$K\mathbf{LeibA}$ associating to each Leibniz algebra 
$\big(\mathfrak{h}, [~,~]\big)$ the augmented Leibniz algebra
$\big(\mathfrak{h}, p, \bar{\mathfrak{h}}, [~,~]_{\bar{\mathfrak{h}}},
\mathrm{ad}'\big)$ where $p:\mathfrak{h}\to \bar{\mathfrak{h}}$ is the canonical
projection and the representation $\mathrm{ad}'$ of the Lie algebra
$\bar{\mathfrak{h}}$ on the Leibniz algebra $\mathfrak{h}$ is defined by
(for all $x,y\in\mathfrak{h}$)
\begin{equation}
    \mathrm{ad}'_{p(x)}(y):= \mathrm{ad}_x(y)=[x,y].
\end{equation}

\subsection{(Augmented) Lie racks}
  \label{SubSecLie racks}

We now restrict to $K=\R$.  
Recall that a {\it pointed manifold} is a pair $(M,e)$ where $M$ is a 
differentiable
manifold and $e$ is a fixed element of $M$. Morphisms of pointed manifold are
base point preserving smooth maps.

Recall that a \emph{Lie rack} is a pointed manifold 
$(M,e)$ equipped with a smooth map $\mathbf{m}:M\times M\to M$ of pointed 
manifolds (i.e.~$\mathbf{m}(e,e)=e$) such that $\mathbf{m}(x,-):M\to M$ 
is a diffeomorphism for all $x\in M$ and 
satisfying the following identities for all $x,y,z\in M$
where the standard notation is $\mathbf{m}(x,y)=x\triangleright y$
\begin{eqnarray}
     e \triangleright x & = & x, \label{EqDefRackETriXEqualsX}\\
     x \triangleright e & = & e, \label{EqDefRackXTriEEqualsE} \\
     x \triangleright ( y \triangleright z) & = &
       (x\triangleright y) \triangleright (x\triangleright z)
         \label{EqSelfDistributivity}
\end{eqnarray}
The last condition (\ref{EqSelfDistributivity}) is called the 
\emph{self distributivity condition}. A morphims of Lie racks
$\phi:$ $(M,e,\mathbf{m})$ $\to$ $(M',e',\mathbf{m}')$ is a map of pointed 
manifolds satisfying for all $x,y\in M$ the condition 
$\phi(x\triangleright y)= \phi(x)\triangleright' \phi(y)$. The class of all 
Lie racks
forms a category called $\mathbf{LieRack}$.
Note that every pointed differentiable manifold $(M,e)$ carries
a \emph{trivial Lie rack structure} defined for all $x,y\in M$ by
\begin{equation}
   x\triangleright_0 y := y,
\end{equation}
and this assignment is functorial.\\
Moreover, any Lie group $G$ becomes a Lie rack upon setting
for all $g,g'\in G$
\begin{equation}
     g\triangleright g' := gg'g^{-1},
\end{equation}
again defining a functor from the category of Lie groups to the 
category of all Lie racks. Examples of racks which are not the conjugation rack underlying a group
abound: Firstly, every conjugation class and every union of conjugation 
classes in a group (defining an immersed submanifold) in a Lie group is 
a Lie rack. 
Then, any Lie rack $(M,e,\triangleright)$ can be \emph{gauged} by any
smooth map $f:(M,e)\to (M,e)$ of pointed manifolds satisfying
for all $x,y\in M$
\[
       f(x\triangleright y) = x\triangleright f(y).
\]
A straight-forward computation shows that the pointed manifold $(M,e)$
equipped with the \emph{gauged} multiplication
$\triangleright_f$ defined by
\[
        x\triangleright_f y := f(x) \triangleright y
\]
is a Lie rack $(M,e,\triangleright_f)$. We refer for more exotic 
examples to \cite{FR92}. 
The following relation to Leibniz algebras is due to M. Kinyon \cite{Kin07}:
\begin{prop}\label{PTangentSpaceLieRackIsLeibnizAlgebra}
  Let $(M,e,\mathbf{m})$ be a Lie rack and $\mathfrak{h}=T_eM$. Define
  the following bracket $[~,~]$ on $\mathfrak{h}$ by
  \begin{equation}\label{EqCompLeibnizBracketFromLieRacks}
      [x,y] = \left.\frac{\partial}{\partial t}T_e L_{a(t)}(y)\right|_{t=0}
  \end{equation}
  where $t\mapsto a(t)$ is any smooth curve defined on an open real interval 
  containing $0$ satisfying $a(0)=e$ and $(da/dt)(0)=x\in\mathfrak{h}$.
  Then we have the following
  \begin{enumerate}
   \item $\big(\mathfrak{h},[~,~]\big)$ is a real Leibniz algebra.
   \item Let $\phi:(M,e,\mathbf{m})\to (M',e',\mathbf{m}')$ be a morphism of
   Lie racks. Then $T_e\phi:\mathfrak{h}\to \mathfrak{h}'$ is a morphism of
   Leibniz algebras.
  \end{enumerate}  
\end{prop}
\begin{prooof}
  Since for each $a\in M$ we have $L_a(e)=e$ it follows that the tangent map 
  $T_eL_a$ maps the tangent space $T_eM$ to $T_eM$ whence the curve
  $t\mapsto T_eL_{a(t)}$ is a curve of $\mathbb{R}$-linear maps $T_eM\to T_eM$
  whence eqn (\ref{EqCompLeibnizBracketFromLieRacks}) defines a well-defined
  real bilinear map $\mathfrak{h}\times \mathfrak{h}\to \mathfrak{h}$.\\
  \textbf{1.}
  Let $x,y,z\in \mathfrak{h}$, and let $t\mapsto a(t)$ and $t\mapsto b(t)$
  two  smooth curves of an open interval (containing $0$) into $M$ such that
  $a(0)=e=b(0)$ and $(da/dt)(0)=x$, $(db/dt)(0)=y$.
  We compute
  \begin{eqnarray*}
    \lefteqn{\big[x,[y,z]\big]} \\ 
      & = &
     \left.\frac{\partial^2}{\partial s\partial t}
       \left(T_e L_{a(s)}\left(T_e L_{b(t)}(z)\right)\right)\right|_{s,t=0}
        =  \left.\frac{\partial^2}{\partial s\partial t}
           T_e\left(L_{a(s)}\circ L_{b(t)} \right)(z)\right|_{s,t=0}\\
       & \stackrel{(\ref{EqSelfDistributivity})}{=}&  
       \left.\frac{\partial^2}{\partial s\partial t}
         T_e\left( L_{a(s)\triangleright b(t)}\circ L_{a(s)}\right)(z)
              \right|_{s,t=0}
        = \left.\frac{\partial^2}{\partial s\partial t}
        \left(  T_eL_{a(s)\triangleright b(t)}
           \left(T_e L_{a(s)}(z)\right)\right)\right|_{s,t=0}\\
      & = & \left.\frac{\partial^2}{\partial s\partial t}
            T_eL_{a(s)\triangleright b(t)}\right|_{s,t=0}
              \big(T_e L_{a(0)}(z)\big) \\
     & &  ~~~+ \left.\frac{\partial}{\partial t}
                    T_eL_{a(0)\triangleright b(t)}\right|_{t=0}
                 \left.\left(\frac{\partial}{\partial s}
                     \left(T_e L_{a(s)}(z)\right)\right)\right|_{s=0}.
  \end{eqnarray*}
  Since $a(0)=e$ we have $T_e L_{a(0)}(z)=z$ and 
  $a(0)\triangleright b(t)=b(t)$ whence the last term equals 
  $\big[y,[x,z]\big]$. Since for each $s$ the curve 
  $t\mapsto a(s)\triangleright b(t)$ is equal to $e$ at $t=0$ we get
  \begin{eqnarray*}
    \left.\frac{\partial^2}{\partial s\partial t}
            T_eL_{a(s)\triangleright b(t)}\right|_{s,t=0}
              \big(T_e L_{a(0)}(z)\big) & = &
        \left[\left. \frac{\partial}{\partial s}T_eL_{a(s)}(y)\right|_{s=0},
               z\right] =\big[[x,y],z\big]
  \end{eqnarray*}
  proving the Leibniz identity.\\
  \textbf{2.} Since $\phi$ maps $e$ to $e'$ its tangent map $T_e\phi$ maps
  $T_eM$ to $T_{e'}M'$. We get for all $x,y\in\mathfrak{h}=T_eM$ where
  $t\mapsto a(t)$ is a smooth curve in $M$ with $a(0)=e$ and $(da/dt)(0)=x$:
  \begin{eqnarray*}
     T_e\phi\big([x,y]\big)
       & = & T_e\phi \left(
         \left.\frac{\partial}{\partial t}T_eL_{a(t)}(y)\right|_{t=0}\right)
       = \left.\frac{\partial}{\partial t}\Big(
            T_e\big(\phi\circ L_{a(t)}\big)(y)\Big)\right|_{t=0} \\
      & = &  \left.\frac{\partial}{\partial t}\Big(
             T_e\big(L'_{\phi(a(t))}\circ \phi\big)(y)\Big)\right|_{t=0} 
        =  \left.\frac{\partial}{\partial t}
           T_{e'} L'_{\phi(a(t))}\right|_{t=0}\big(T_e\phi(y)\big) \\
     & = & \big[T_e\phi(x),T_e\phi(y)\big].
  \end{eqnarray*}
\end{prooof}
Let $\mathbb{R}\mathbf{Leib}_{fd}$ denote the category of all finite-dimensional
real Leibniz algebras.
The preceding proposition shows that there is a functor
$T_*\mathcal{R}:\mathbf{LieRack}\to \mathbb{R}\mathbf{Leib}_{fd}$ which 
associates to any Lie rack $(M,e,\triangleright)$ its tangent space
$T_*\mathcal{R}(M):=T_eM$ at the distinguished point $e\in M$ 
equipped with the Leibniz bracket
eqn (\ref{EqCompLeibnizBracketFromLieRacks}).

Furthermore, recall that an \emph{augmented Lie rack} (see \cite{FR92})
$(M,\phi,G,\ell)$ consists of a pointed differentiable manifold 
$(M,e_M)$,
of a Lie group $G$, of a smooth map $\phi:M\to G$ (of pointed manifolds),
and of a
smooth left $G$-action $\ell:G\times M\to M$ 
(written $(g,x)\mapsto \ell(g,x)=\ell_g(x)=gx$) such that for all
$g\in G$, $x\in M$
\begin{eqnarray}
     ge_M & = & e_M, \label{EqDefAugRackGActionPreservesE} \\
     \phi(gx) & = & g\phi(x)g^{-1}. 
     \label{EqDefAugRackPhiIntertwinesGActionAndConjugation}
\end{eqnarray}
It is a routine check that the multiplication $\triangleright$ on $M$
defined for all $x,y\in M$ by
\begin{equation}
  x\triangleright y  :=  \ell_{\phi(x)}(y)
\end{equation}
satisfies all the axioms (\ref{EqDefRackETriXEqualsX}), 
(\ref{EqDefRackXTriEEqualsE}), and (\ref{EqSelfDistributivity})
of a Lie rack, thus making $(M,e_M,\triangleright)$ into a 
Lie rack such that the map $\phi$ is a morphism of Lie racks, i.e.
for all $x,y\in M$
\begin{equation}
     \phi(x\triangleright y)= \phi(x)\phi(y)\phi(x)^{-1}.
\end{equation}
A morphism $(\Psi,\psi):(M,\phi,G,\ell)\to (M',\phi',G',\ell')$ of 
augmented
Lie racks is a pair of maps of pointed differentiable manifolds
$\Psi:M\to M'$ and $\psi:G\to G'$ such that $\psi$ is homomorphism
 of Lie groups and such that all reasonable diagrams commute, viz:
 for all $g\in G$
\begin{eqnarray}
   \phi'\circ \Psi & = & \psi\circ \phi, \\
   \Psi\circ \ell_{g} & = & \ell'_{\psi(g)}\circ \Psi.
\end{eqnarray}
The class of all augmented Lie racks thus forms a category $\mathbf{LieRackA}$
with the obvious forgetful functor $F:\mathbf{LieRackA}\to \mathbf{LieRack}$.
Note that the trivial Lie rack structure of a pointed manifold $(M,e)$
comes from an augmented Lie rack over the trivial Lie group $G=\{e\}$.

\section{Dirty integration of Leibniz algebras}

\subsection{The main theorem}

Let $\big(\mathfrak{h},p,\mathfrak{g},[~,~]_\mathfrak{g},\dot{\rho}\big)$
be an augmented Leibniz algebra, let $\mathfrak{g}'$ be the Lie ideal
$p(\mathfrak{h})$ of $\mathfrak{g}$, and let $\mathfrak{z}:=\mathrm{Ker}(p)$
which --we recall-- is a two-sided ideal of the Leibniz algebra $\mathfrak{h}$ 
lying in the left centre of $\mathfrak{h}$. Let furthermore $G$ (resp. $G'$)
be a connected simply connected Lie group whose Lie algebra is isomorphic
to $\mathfrak{g}$ (resp.~$\mathfrak{g}'$). Since $G$ is connected and 
simply connected,
its adjoint represention $\mathrm{Ad}_G$ preserves the ideal $\mathfrak{g}'$
of its Lie algebra $\mathfrak{g}$, whence there is a Lie group homomorphism
$g\mapsto \mathrm{A}'_g$ of $G$ into $\mathrm{Aut}_0(\mathfrak{g}')$, 
the component of the identity of
the Lie group of all automorphisms of the Lie algebra $\mathfrak{g}'$.
Since this latter Lie group is well known to be isomorphic to
$\mathrm{Aut}_0(G')$, the connected
component of the identity of the topological group of all Lie group automorphisms
of $G'$ (which also is a Lie group), see e.g.~\cite{Hoc52} for details,
there is a unique Lie group homomorphism
\[
     \mathrm{I}':G\to \mathrm{Aut}_0(G'):g\mapsto 
     \big(g'\mapsto \mathrm{I}'_g(g')\big)
\]
such that $T_{e}\mathrm{I}'_g = \mathrm{A}_g$ for all $g\in G$.
Moreover, the
injection
$\mathfrak{g}'\to \mathfrak{g}$ induces a unique immersive Lie group homomorphism
$\iota:G'\to G$ whose image is an analytic
normal subgroup of $G$ whence
$G$ acts on $\iota(G')$ by conjugations, and we have
for all $g\in G$ and all $g'\in G'$
\[
      (\iota\circ \mathrm{I}'_g )(g')= g\iota(g')g^{-1}.
\]
Next, let $\rho:G\times \mathfrak{h}\to \mathfrak{h}$ be the unique
representation of $G$ on $\mathfrak{h}$ such that for all $\xi\in\mathfrak{g}$
and $x\in\mathfrak{h}$
\begin{equation}
     \left.\frac{d}{dt}\rho_{\exp(t\xi)}(x)\right|_{t=0} =\dot{\rho}_\xi(x).
\end{equation}
We get for all $g\in G$ and $x\in\mathfrak{h}$:
\begin{equation}
     p\big(\rho_g(x)\big) = \mathrm{A}'_g\big(p(x)\big).
\end{equation}

\noindent 
The main theorem of this article reads:

\begin{theorem}  \label{main_theorem}
  With the above hypotheses and notations we have the following:
  \begin{enumerate}
  \item There is a smooth map $\mathsf{s}:G'\to \mathfrak{g}'$ having the 
  following 
    properties:
    \begin{eqnarray}
       \mathsf{s}(e) & = & 0, \\
       T_e\mathsf{s} & = & \mathrm{id}_{\mathfrak{g}'}, \\
       \forall~g\in G,~\forall~g'\in G':~~~
       \mathsf{s}\big(\mathrm{I}'_g(g')\big) & = & 
         \mathrm{A}'_{g}\big(\mathsf{s}(g')\big).
    \end{eqnarray}

  \item Consider $p:\mathfrak{h}\to \mathfrak{g}'$ as a fibre bundle over
       $\mathfrak{g}'$ (it is an \textbf{affine} bundle 
       with typical fibre $\mathfrak{z}$ over $\mathfrak{g}'$), and form the 
       pulled-back fibre bundle
       \begin{equation}
            M := \mathsf{s}^*\mathfrak{h}
              = \big\{(x,g')\in \mathfrak{h}\times G'~|~
                       p(x)=\mathsf{s}(g')\big\}
                       \stackrel{\phi}{\rightarrow} G'
       \end{equation}
       over $G'$ having $(0,e_{G'})$ as a distinguished point.
       There is a canonical $G$-action $\ell$ on $M$ induced by 
       $\rho$ on $\mathfrak{h}$ and by
       $\mathrm{I}'$ on $G'$ such that
       $\big((M,(0,e_{G'}),\iota\circ \phi, G,\ell\big)$ is an augmented 
       Lie rack.
       
  \item The induced Leibniz algebra structure on the tangent space
  $T_{(0,e_{G'})}M$ is isomorphic to $\big(\mathfrak{h},[~,~]_\mathfrak{h}\big)$.
  
  \item In the particular case $\mathfrak{g}=\mathfrak{g}'$ and $G=G'$ the above
   construction gives a surjective projection $M\to G=G'$. If furthermore
   $\mathfrak{z}=\{0\}$ the above construction reduces to 
   the usual conjugation Lie rack on $G=G'$.
   \end{enumerate}
\end{theorem}

\subsection{Proof of the main theorem}

\textbf{1.} 
According to Proposition \ref{PAdInvariantChartDomainExpAndBumpFunction} 
(which we separately show further down in the Appendix)
there are two open neighbourhoods
$\mathcal{U}'_{(\pi/2)\mathbf{i}}\subset\mathcal{U}'_{\pi\mathbf{i}}$ 
of $0\in\mathfrak{g}'$ which are both
$\mathrm{Aut}_0(\mathfrak{g}')$-invariant and on which the restriction of the 
exponential map 
$\exp_{G'}$ is a diffeomorphism onto the $\mathrm{Aut}_0(G')$-invariant open
neighbourhoods 
$\mathcal{V}'_{(\pi/2)\mathbf{i}}\subset\mathcal{V}'_{\pi\mathbf{i}}$ of the 
unit element $e'$ of $G'$.
Moreover, again by 
Proposition \ref{PAdInvariantChartDomainExpAndBumpFunction} there is an 
$\mathrm{Aut}_0(\mathfrak{g}')$-invariant bump function 
$\gamma':\mathfrak{g}'\to [0,1]$ whose support is in 
$\mathcal{U}'_{\pi\mathbf{i}}$ and which is equal to $1$ on 
$\mathcal{U}'_{(\pi/2)\mathbf{i}}$. Let us define the following map
$\mathsf{s}:G'\to \mathfrak{g}'$ by
\begin{equation}
   \mathsf{s}(g'):= \left\{
      \begin{array}{cl}
         \gamma'\big(\exp_{G'}^{-1}(g')\big)\exp_{G'}^{-1}(g')
            & \forall~g'\in \mathcal{V}'_{\pi\mathbf{i}} \\
            0 & \forall~g'\not \in \mathcal{V}'_{\pi\mathbf{i}}
      \end{array}\right.
\end{equation}
It is clear that $\mathsf{s}$ is a well-defined smooth map $G'\to \mathfrak{g}'$.
Moreover $\mathsf{s}(e_{G'})=0$ by the properties of the exponential map,
and (setting $e=e_{G'}$) for all $\zeta\in\mathfrak{g}'=T_{e}G'$
\begin{eqnarray*}
  T_e\mathsf{s} (\zeta) 
    & = &\frac{d}{ds}\left(
     \gamma'\big(\exp_{G'}^{-1}(\exp_{G'}(s\zeta))\big)
     \exp_{G'}^{-1}(\exp_{G'}(s\zeta))\right)|_{s=0}\\
   & = & \frac{d}{ds}\left(s\zeta\right)|_{s=0} = \zeta
\end{eqnarray*}
because the bump function $\gamma'$ is constant equal to $1$ near $0$.
Hence $T_e\mathsf{s}=\mathrm{id}_{\mathfrak{g}'}$. Finally, since
for each $g\in G$ and $\zeta\in\mathfrak{g}'$
\[
    \mathrm{I}'_g\big(\exp_{G'}(\zeta)\big)
       = \exp_{G'}\big(\mathrm{A}'_g(\zeta)\big)
\]
and since $\gamma'$ is invariant under the action $g\mapsto\mathrm{A}'_g$ of
$G$ on $\mathfrak{g}'$
we get for all $g\in G$ and $g'\in G'$
\[
   \mathsf{s}\big(\mathrm{I}'_g(g')\big)= \mathrm{A}'_g\big(\mathsf{s}(g')\big),
\]
proving the first statement of the theorem.

\textbf{2.} Since $p:\mathfrak{h}\to \mathfrak{g}'$ is a surjective linear map,
it is a surjective submersion whose fibre over $0\in\mathfrak{g}'$ is equal to 
$\mathrm{Ker}(p)=\mathfrak{z}$, and whose fibre over any $\zeta\in\mathfrak{g}'$
is the affine subspace $p^{-1}(\{\zeta\})$ of $\mathfrak{h}$. Choosing any
vector space complement $\mathfrak{b}$ to $\mathfrak{z}$ in $\mathfrak{h}$
leads to differential geometric trivialization over the global chart domain 
$\mathfrak{g}'$ of $\mathfrak{g}'$.
Hence $p:\mathfrak{h}\to \mathfrak{g}'$ is a fibre bundle over $\mathfrak{g}'$
with typical fibre $\mathfrak{z}$, and therefore the 
pull-back $\mathsf{s}^*\mathfrak{h}=M$ is a well-defined fibre bundle over $G'$
with typical fibre $\mathfrak{z}$. Recall that the projection $\phi:M\to G'$
is given by the restriction of the projection on the second factor 
$\mathfrak{h}\times G'\to G'$ to the submanifold 
$M\subset \mathfrak{h}\times G'$.\\
Since $\mathsf{s}(e_{G'})=0=p(0)$ it follows that the point $(0,e_{G'})$
is in $M$, and clearly $\phi(0,e_{G'})=e_{G'}$.\\
There is a canonical diagonal $G$-action $\hat{\ell}$ on 
$\mathfrak{h}\times G'$ defined by 
$\hat{\ell}_g(x,g')= \big(\rho_g(x),\mathrm{I}'_g(g')\big)$. As for any
$(x,g')\in M$ we have by definition $p(x)=\mathsf{s}(g')$, we get
\[
    \mathsf{s}\big(\mathrm{I}'_g(g')\big) = 
      \mathrm{A}'_g\big(\mathsf{s}(g')\big) = 
      \mathrm{A}'_g\big(p(x)\big)=p\big(\rho_g(x)\big)
\]
proving that for any $(x,g')\in M$ the point $\hat{\ell}_g(x,g')\in M$,
whence $\hat{\ell}$ restricts to a well-defined $G$-action $\ell$ on $M$.
Clearly $\ell_g(0,e_{G'})=(0,e_{G'})$. Moreover, for any $(x,g')\in M$
and $g\in G$ we have
\begin{eqnarray*}
 (\iota\circ\phi)\big(\ell_g(x,g')\big)& = &
 \iota\Big(\phi\big(\rho_g(x),\mathrm{I}'_g(g')\big)\Big)=
   \iota\big(\mathrm{I}'_g(g')\big) \\
   &= &g\iota(g')g^{-1}
   =g\big((\iota\circ \phi)(x,g')\big)g^{-1},
\end{eqnarray*}
showing that $\big((M,(0,e_{G'}),\iota\circ \phi, G,\ell\big)$ is an augmented 
       Lie rack.
       
\textbf{3.} According to the definition of the pull-back, the tangent space
of $M$ at $(e_{G'},0)$ is given by all the pairs $(x,\zeta)\in
\mathfrak{h}\times\mathfrak{g}'  = T_0\mathfrak{h}\times T_{e_{G'}}G'$ such that
\[
   T_{e_{G'}}\mathsf{s}(\zeta)=p(x),~~\mathrm{hence}~
        \zeta=p(x)
\]
because of $T_{e_{G'}}\mathsf{s}=\mathrm{id}_{\mathfrak{g}'}$. It follows that
the linear map 
$\theta_\mathfrak{h}:\mathfrak{h}\to \mathfrak{h}\times \mathfrak{g}'$ given by
\[
   \theta_\mathfrak{h}(x)= \big(x,p(x)\big)
\]
is an isomorphism of the vector space $\mathfrak{h}$
onto the tangent space $T_{(e_{G'},0)}M$. We get for all $g\in G$ and for all 
$y\in\mathfrak{h}$
\begin{eqnarray*}
  T_{(e_{G'},0)}\ell_g\big(y,p(y)\big)  
  & = & \left.\frac{d}{dt}\Big(\rho_g(ty),
  \mathrm{I}'_g\big(\exp\big(t p(y)\big)\big)\Big)\right|_{t=0}
  =\Big(\rho_g(y),\mathrm{A}'_g\big(p(y)\big)\Big).
\end{eqnarray*}
Now for all $x,y\in\mathfrak{h}$, we get for the Leibniz bracket on
the tangent space $T_{(e_{G'},0)}M$ of the augmented Lie rack
$\big((M,(0,e_{G'}),\iota\circ \phi, G,\ell\big)$
\begin{eqnarray*}
  \big[\theta_\mathfrak{h}(x),\theta_\mathfrak{h}(y)\big] & =&
  \Big[\big((x,p(x)\big),\big(y,p(y)\big)\Big] \\
  & = &\left.\frac{d}{dt} T_{(e_{G'},0)}\ell_{\exp\big(tp(x)\big)}
                \big(y,p(y)\big)\right|_{t=0} \\
   & = &\left.\frac{d}{dt}
   \Big(\rho_{\exp\big(tp(x)\big)}(y),
     \mathrm{A}'_{\exp\big(tp(x)\big)}\big(p(y)\big)\Big)\right|_{t=0} \\
     & = & \Big( [x,y]_\mathfrak{h},
              \big[p(x),p(y)\big]_{\mathfrak{g}'}\Big)
        = \Big( [x,y]_\mathfrak{h},
              p\big([x,y]_\mathfrak{h}\big)\Big)
        = \theta_\mathfrak{h}\big([x,y]_\mathfrak{h}\big).
\end{eqnarray*}
showing that the induced Leibniz structure form the augmented Lie rack
is isomorphic with the Leibniz bracket $[~,~]_\mathfrak{h}$ on
$\mathfrak{h}$.

\noindent \textbf{4.} This is immediate.\hfill $\Box$

\appendix

\section{Automorphism invariant chart domains for the exponential map 
         of a Lie group}
         
\subsection{Zeros of Polynomials}

All the results in this section are classical, and the methods
had been inspired e.g. by the article \cite{NP94}.

Let $\mathbb{C}[\lambda]^1_n$ denote the space of all 
\emph{monic complex polynomials of degree $n$}. Since every
such polynomial $f$ is of the general form 
\begin{equation}\label{EqDefPolynomialMonic}
   f(\lambda) =\lambda^n +\sum_{r=1}^n a_r \lambda^{n-r}
\end{equation}
with $a:=(a_1,\ldots,a_n)\in\mathbb{C}^n$ it is clear that
$\mathbb{C}[\lambda]^1_n$ is an affine space with associated
complex vector space $\mathbb{C}^n$, whence $\mathbb{C}[\lambda]^1_n$
is homeomorphic to $\mathbb {C}^n$. We put the norm
$||f||=||a||:=\sum_{j=1}^n |a_j|$ on it. Consider now the map 
$\mathcal{T}:\mathbb{C}^n\to \mathbb{C}[\lambda]^1_n$ given by
\begin{equation}
   \mathcal{T}(z)(\lambda):= (\lambda-z_1)\cdots (\lambda -z_n)
\end{equation}
so that the zeros of $\mathcal{T}(z)$ are given by 
$z_1,\ldots,z_n\in \mathbb{C}$ (where there can be repetitions, i.e.
multiple roots). There is the well-known classical formula expressing
the coefficients of $\mathcal{T}(z)$ in terms of Newton's elementary
symmetric polynomials, i.e.
\begin{equation}
   \mathcal{T}(z)(\lambda)=\lambda^n
       +\sum_{r=1}^n(-1)^r
       \big(\sum_{1\leq i_1<\cdots<i_r\leq n}z_{i_1}\cdots z_{i_r}\big)
          \lambda^{n-r}.
\end{equation}
It follows that $\mathcal{T}$ is a complex analytic map, hence continuous.
Moreover the Fundamental Theorem of Algebra states that $\mathcal{T}$ is 
surjective, and elementary algebra of polynomials shows that
\begin{equation}\label{EqCompInjectivityFailureOfMathcalT}
   \mathcal{T}(z)=\mathcal{T}(z')~~~\mathrm{if~and~only~if}~~~
      \exists~\sigma\in S_n~\mathrm{such~that~}z'= z.\sigma
\end{equation}
where a permutation $\sigma$ in the symmetric group $S_n$ acts from the right
on $\mathbb{C}^n$ in the usual way, i.e. $(z_1,\ldots,z_n).\sigma
= (z_{\sigma(1)},\ldots,z_{\sigma(n)})$. Let $\mathbb{C}^n/S_n$ be the 
space of all $S_n$-orbits equipped with the quotient topology, and let
$\pi_n:\mathbb{C}^n\to\mathbb{C}^n/S_n$ be the canonical projection. 
Eqn (\ref{EqCompInjectivityFailureOfMathcalT}) implies that the map 
$\mathcal{T}$ descends to a well-defined continuous map 
$T:\mathbb{C}^n/S_n\to \mathbb{C}[\lambda]^1_n$ which is bijective.
It is classical, but a bit less well-known that \emph{$T$ is a homeomorphism}:
\begin{prop}
   With the above notations we have the following:
   \begin{enumerate}
    \item The map $\mathcal{T}$ is a closed continuous map.
    \item The map $T$ is a homeomorphism, and $\mathcal{T}$ is also an open
       map.
   \end{enumerate}
\end{prop}
\begin{prooof}
  We need first the following elementary estimate: Consider a monic polynomial
  $f\in \mathbb{C}[\lambda]^1_n$ in the form (\ref{EqDefPolynomialMonic}),
  and let $\mu$ be a root of $f$, then in case $|\mu|\geq 1$ we get
  from the equation $f(\mu)=0$ 
  \[
     \mu^n= -\sum_{r=1}^na_r\mu^{n-r}~~\mathrm{hence}~~~
     |\mu|^n \leq \sum_{r=1}^n|a_r||\mu|^{n-r}~~\mathrm{hence}~~~
     |\mu| \leq ||a||
  \]
  where we have multiplied the second term by $|\mu|^{1-n}>0$ 
  and used the fact that $|\mu|^{1-r}\leq 1$ for all $1\leq r\leq n$.
  This estimate implies the weaker estimate
  \begin{equation}\label{InEqZerosInTermsOfCoefficients}
      |\mu| \leq \mathrm{max}\{1, ||a|| \}
  \end{equation}
  which clearly also holds for the other case $|\mu|\leq 1$.\\
 \textbf{1.} Let $F$ be a closed subset of $\mathbb{C}^n$, and consider a 
   sequence
   $(f_k)_{k\in\mathbb{N}}$ in $\mathcal{T}(F)$ converging to a monic polynomial
   $f\in \mathbb{C}[\lambda]^1_n$ (where we use the above norm $||~||$ on the 
   coefficients to define the convergence). We have to show that there
   is a $z\in F$ such that $f=\mathcal{T}(z)$. First, since each 
   $f_k\in \mathcal{T}(F)$, there is a sequence $(z^{(k)})_{k\in\mathbb{N}}$
   of elements of $F$ such that for all $k\in\mathbb{N}$, we have
   $f_k=\mathcal{T}(z^{(k)})$. Denote by $a^{(k)}\in\mathbb{C}^n$ the 
   coeffients of the polynomial $f_k$, and let $a\in\mathbb{C}^n$ be the vector 
   of coeffients
   of the polynomial $f$. By hypothesis $a^{(k)}\to a$ when $k\to\infty$,
   hence the sequence of norms $||a^{(k)}||$ converges to $||a||$ and is 
   therefore bounded. But the above estimate 
   (\ref{InEqZerosInTermsOfCoefficients}) implies that the sequence
   $(z^{(k)})_{k\in\mathbb{N}}$ of zero vectors of $f_k$ is a bounded subset
   of $\mathbb{C}^n$. By the Bolzano-Weierstrass Theorem, there is a 
   subsequence $\big(z^{(k_l)}\big)_{l\in\mathbb{N}}$ of the above sequence
   which converges to $z\in\mathbb{C}^n$. Since all the vectors $z^{(k_l)}$
   are in the closed set $F$, it follows that the limit $z$ lies also in $F$.
   By the continuity of $\mathcal{T}$ it follows that
   \[
      f_{k_l}=\mathcal{T}\big(z^{(k_l)}\big)\to \mathcal{T}(z)~~~(l\to\infty),
   \]
   and since the subsequence $(f_{k_l})_{l\in\mathbb{N}}$ converges to
   $f$, it follows by the uniqueness of limits that $f=\mathcal{T}(z)$,
   and $\mathcal{T}$ is a closed map.\\
   \textbf{2.} We shall show that $T$ is a closed map which will imply
   that its inverse
   is continuous: Let $F$ be a closed subset of $\mathbb{C}^n/S_n$. Then
   $\pi_n\big(\pi_n^{-1}(F')\big)=F'$ because $\pi_n$ is surjective, and
   \[
       T(F')=T\big(\pi_n\big(\pi_n^{-1}(F')\big)\big)=
         \mathcal{T}\big(\pi_n^{-1}(F')\big)
   \]
  is a closed subset of $\mathbb{C}[\lambda]^1_n$, because $\pi_n^{-1}(F')$
  is a closed subset of $\mathbb{C}^n$ thanks to the continuity of $\pi_n$,
  and because $\mathcal{T}$ is a closed map. In order to prove the second half,
  observe that $\pi_n$ is an open map: If $U\subset\mathbb{C}^n$ is open, then
  $\pi_n^{-1}\big(\pi_n(U)\big)=\bigcup_{\sigma\in S_n} U.\sigma$ is a union
  of the open sets $U.\sigma$ of $\mathbb{C}^n$, and therefore open whence
  $\pi_n(U)$ is open by definition of the quotient topology. Since $T$ is a
  homeomorphism, it is an open map (its inverse is continuous), and therefore
  $\mathcal{T}=T\circ \pi_n$ is open as a composition of open maps. 
\end{prooof}

\noindent The following corollary will be important in the proof of
Theorem \ref{main_theorem}:
\begin{cor}\label{CPolynomialsHavingAllRootsInOpenSetsAreOpen}
 Let $U$ be an open and $F$ be a closed subset of $\mathbb{C}$. Then for each 
 positive integer $n$, the set of all those monic polynomials of degree $n$
 having all their roots in $U$ (resp.~in $F$) is an open (resp.~closed) subset
 of $\mathbb{C}[\lambda]_n^1$.
\end{cor}
\begin{prooof}
 Apply the map $\mathcal{T}$ to the open (resp.~closed) subset 
 $U^{\times n}$ (resp.~$F^{\times n}$) of $\mathbb{C}^n$ and use the fact that
 $\mathcal{T}$ is an open and closed map.
\end{prooof}

\subsection{Injectivity domains of the exponential of linear maps}

Let $E$ be a real vector space of dimension $n$, and let 
$\mathcal{A}=\mathrm{Hom}_\mathbb{R}(E,E)$ be the real vector space of all
$\mathbb{R}$-linear maps $E\to E$. Let $g:\mathbb{C}\to\mathbb{C}$ be a
holomorphic function whose power series expansion around $0$ has
real coefficients, $g(z)=\sum_{r=0}^\infty g_r z^r$, $g_r\in\mathbb{R}$
for all $r\in\mathbb{N}$. Then for any $X\in\mathcal{A}$ the series
$g(X):=\sum_{r=0}^\infty g_r X^r$ is well-known to converge to a
well-defined element in $\mathcal{A}$ (where we write $XY$ for the composition 
of linear maps $X\circ Y$, and set
$X^0:=I:=\mathrm{id}_E$). An important example is
the exponential function
\begin{equation}
   \exp:\mathcal{A} \to \mathcal{A}: X\mapsto \exp(X):= e^X
         :=\sum_{r=0}\frac{1}{r!}X^r,
\end{equation}
but also the series
\begin{equation}\label{EqDefFunctionHDerivativeOfExponential}
   h(X) : = \frac{I-e^{-X}}{X}:= \sum_{r=0}^\infty \frac{(-1)^r}{(r+1)!}X^r,
\end{equation}
related to the derivative of the exponential function.
For any strictly positive real number $\tau$, let $S_{\tau\mathbf{i}}$ be the open 
 strip
\begin{equation}\label{EqDefOpenStripOfC}
   S_{\tau\mathbf{i}}:=\{z=\alpha+\mathbf{i}\beta\in\mathbb{C}~|~|\beta|<\tau\},
\end{equation}
and let
\begin{equation}
   \mathcal{A}_{\tau\mathbf{i}}:=
      \{X\in\mathcal{A}~|~\mathrm{all~the~eigenvalues~of~}X\mathrm{~are~in~}
           S_{\tau\mathbf{i}}\}.
\end{equation}
\begin{prop} \label{PInjectivityRegularityExpLinearMaps}
For any positive integer $n$, we have the following:
  \begin{enumerate}
  \item For any positive real number $\tau>0$, the subset
    $\mathcal{A}_{\tau\mathbf{i}}$ is an open subset of $\mathcal{A}$.
 \item For all $X\in \mathcal{A}_{2\pi\mathbf{i}}$, the linear map 
   $h(X)$ is invertible. 
   \item The restriction of the exponential map 
 to the open subset $\mathcal{A}_{\pi\mathbf{i}}$ of $\mathcal{A}$ is injective.
  \end{enumerate}
\end{prop}
\begin{prooof}
 For any $X\in\mathcal{A}$, recall the Jordan decomposition 
 \[
       X= X_S+X_N
 \]
into its semisimple part $X_S\in\mathcal{A}$ and its nilpotent part 
$X_N\in\mathcal{A}$ where of course
$X_S$ is diagonalizable over $\mathbb{C}$ and $X_N$ is nilpotent, and
$X_S$ and $X_N$ are polynomials of $X$, hence commute with each other.
The Jordan decomposition is well-known to be unique. Now take any
holomorphic function $g:\mathbb{C}\to\mathbb{C}$ whose power series expansion 
around $0$ has real coefficients. It follows that the function of two complex 
variables 
\[
      (z,w)\mapsto \frac{g(z+w)-g(z)}{w}
\]
is a well-defined holomorphic function $\mathbb{C}^2\to\mathbb{C}$
(with $(z,0)$ sent to $f'(z)$), whence
the decomposition 
\[
     g(X)=g\big(X_S+X_N\big) =g\big(X_S\big) + 
            X_N\frac{g\big(X_S+X_N\big)-g\big(X_S\big)}{X_N}
\]
is well-defined. If $\mu_1,\ldots,\mu_n\in\mathbb{C}$ are the eigenvalues of
$X$ hence of $X_S$, then $g(\mu_1),\ldots,g(\mu_n)$ are the eigenvalues
of $g\big(X_S\big)$ (with the same eigenvectors) whence $g\big(X_S\big)$
is obviously semisimple, and since
the second
summand is equal to the product of the nilpotent map $X_N$ and another map 
commuting with $X_N$, it follows that the preceding equation gives us
the Jordan decomposition of $g(X)$: 
\begin{equation}\label{EqCompSemisimpleNilpotentPartOfGOfX}
   g(X)_S= g\big(X_S\big),~~~\mathrm{and}~~~
      g(X)_N= X_N\frac{g\big(X_S+X_N\big)-g\big(X_S\big)}{X_N}
\end{equation}
\textbf{1.} Consider the map $\chi:\mathcal{A}\to \mathbb{C}[\lambda]^1_n$ which sends
each $X\in\mathcal{A}$ to its characteristic polynomial 
$\lambda \mapsto \det(\lambda I - X)$. This map is clearly polynomial in the 
coefficients of $X$, and therefore continuous. By Corollary 
\ref{CPolynomialsHavingAllRootsInOpenSetsAreOpen}, the set of all monic complex 
polynomials all of whose roots are in $S_{\tau\mathbf{i}}$  is the open set 
$\mathcal{T}(S_{\tau\mathbf{i}}^{\times n})$ of $\mathbb{C}[\lambda]^1_n$, and hence
$\mathcal{A}_{\tau\mathbf{i}}=\chi^{-1}
\big(\mathcal{T}(S_{\tau\mathbf{i}}^{\times n})\big)$
is an open subset of $\mathcal{A}$.\\
\textbf{2.} We get that
\[
    \det\big(h(X)\big)=\det\big(h(X)_S\big)=\det\big(h\big(X_S\big)\big).
\]
Let $\mu_1,\ldots,\mu_n$ be the eigenvalues of $X_S$ (repetitions due to
multiplicities may occur). Then
\[
\det\big(h\big(X_S\big)\big)=h(\mu_1)\cdots h(\mu_n).
\]
In case $\mu_r=0$, then $h(\mu_r)=h(0)=1\neq 0$. Let 
$\mu_r=\alpha+\beta\mathbf{i}\neq 0$. Then if $h(\mu_r)=0$, we get
$\mu_rh(\mu_r)=0$, hence
\[
     0=1-e^{-\alpha-\beta\mathbf{i}}
\]
which is the case iff $\alpha_r=0$ and $\beta=2k\pi$ for some integer $k$.
Since by hypothesis $|\beta|<2\pi$, we would necessarily have $\beta_r=0$,
contradiction! Hence all the complex numbers $h(\mu_r)$ are different from
zero, and therefore $h(X)$ is invertible.\\
\textbf{3.}
Let $X,Y\in\mathcal{A}$ such that $e^X=e^Y$. Then we get
\[
     e^{X_S}=\big(e^X\big)_S=\big(e^Y\big)_S=e^{Y_S}.
\]
Let $\mu_1,\ldots,\mu_n\in\mathbb{C}$ be the eigenvalues of $X_S$ and
let $\nu_1,\ldots,\nu_n\in\mathbb{C}$ be the eigenvalues of $Y_S$ (with possible
repetitions due to multiple eigenvalues). Then both
$e^{\mu_1},\ldots,e^{\mu_n}$ and $e^{\nu_1},\ldots,e^{\nu_n}$ are the eigenvalues
of $e^{X_S}=e^{Y_S}$, and we can assume after a possible permutation that
\[
 e^{\mu_1}=e^{\nu_1},~\ldots,~e^{\mu_n}=e^{\nu_n}
\]
Decomposing into real and imaginary part, i.e.
$\mu_s=\alpha_s+\beta_s\mathbf{i}$ and 
$\nu_t=\gamma_t+\delta_t\mathbf{i}$ for all integers $1\leq s,t\leq n$, we can 
conclude that there exist integers $k_1,\ldots,k_n$ such that for all
$1\leq s \leq n$
\[
   \gamma_s=\alpha_s,~~~\mathrm{and}~~~\delta_s=\beta_s + 2\pi k_s.
\]
Hence if both $X$ and $Y$ (and thus $X_S,Y_S$) are in 
$\mathcal{A}_{\mathbf{i}\pi}$, then 
\[
   |\delta_s-\beta_s|\leq |\delta_s|+|\beta_s| < 2\pi
\]
and if $e^{X_S}=e^{Y_S}$, then $X_S$ and $Y_S$ have the same eigenvalues and
the same multiplicities. Let $\hat{\mu}_1,\ldots,\hat{\mu}_k\in\mathbb{C}$
be the $k$ pairwise different eigenvalues of $X_S$ and of $Y_S$. Then we can 
write
\[
    X_S=\sum_{r=1}^k \hat{\mu}_rP_r,~~ Y_S=\sum_{r=1}^k \hat{\mu}_rQ_r
\]
where $P_r,Q_r: E\otimes_\mathbb{R} \mathbb{C}\to E\otimes_\mathbb{R} \mathbb{C}$
are the projections on the eigenspace associated to the eigenvalue 
$\hat{\mu}_r$ (recall that $P_rP_s=\delta_{rs}P_r$ and $Q_rQ_s=\delta_{rs}Q_r$
for all $1\leq r,s \leq k$). We get
\[
    \sum_{r=1}^k e^{\hat{\mu}_r}P_r   =e^{X_S}=e^{Y_S}
         = \sum_{r=1}^k e^{\hat{\mu}_r}Q_r,
\]
and since all the $k$ complex numbers $e^{\hat{\mu}_1},\ldots,e^{\hat{\mu}_k}$
are pairwise different by the fact that $\hat{\mu}_1,\ldots,\hat{\mu}_k\in
S_{\pi\mathbf{i}}$, it follows that for each $1\leq r\leq k$ the 
projection $P_r$ is the unique projection of $e^{X_S}=e^{Y_S}$ on the 
generalized
eigenspace associated to the eigenvalue $e^{\hat{\mu}_r}$ whence $P_r=Q_r$,
and therefore $X_S=Y_S$. Finally we have --according to eqn 
(\ref{EqCompSemisimpleNilpotentPartOfGOfX})--
\begin{eqnarray*}
  e^X(I-e^{-X_N}) & =&  e^X X_N h(X_N)  =  
    X_N\frac{e^{X_S+X_N}-e^{X_S}}{X_N} = \big(e^X\big)_N \\
      & =& \big(e^Y\big)_N = Y_N\frac{e^{Y_S+Y_N}-e^{Y_S}}{Y_N}
       =  e^Y Y_N h(Y_N) = e^Y(I-e^{-Y_N}) 
\end{eqnarray*}
whence --since $e^X=e^Y$-- it follows that $e^{X_N}=e^{Y_N}$.
Since $e^{X_N}$ is a polynomial in $X_N$, and $e^{X_N}-I$ is nilpotent,
we may apply the logarithmic series $\log(1+z)=\sum_{r=0}^\infty (-z)^r/(r+1)$
to $e^{X_N}-I=e^{Y_N}-I$ to get $X_N=Y_N$. It follows that 
\[
     X=X_S+X_N=Y_S+Y_N=Y,
\]
and the restriction of the exponential function to $\mathcal{A}_{\pi\mathbf{i}}$
is injective.
\end{prooof}

\subsection{Automorphism invariant chart domains for the exponential map 
         of a Lie group}

The results and techniques from this section are largely due to \cite{LT66}
and Definition 4.8 and Proposition 4.9 from \cite{DH97}.

Let $\big(\mathfrak{g},[~,~]\big)$ be a finite-dimensional real Lie algebra,
and let $G$ be a connected, simply connected Lie group whose Lie algebra
is isomorphic to $\mathfrak{g}$ (recall that $G$ is unique up to isomorphism
of Lie groups). Moreover, let $\mathrm{Aut}_0(G)$ be the connected component
of the identity of the topological group (w.r.t. the compact-open topology)
of all smooth Lie group automorphisms
of $G$. Then it is known (see e.g.~\cite{Hoc52}) that $\mathrm{Aut}_0(G)$
is a connected Lie group isomorphic to the connected component of the identity,
$\mathrm{Aut}_0(\mathfrak{g})$, of
the Lie group $\mathrm{Aut}(\mathfrak{g})$ of all automorphisms of the Lie 
algebra $\mathfrak{g}$ where the canonical map $\mathrm{Aut}_0(G)\to
\mathrm{Aut}_0(\mathfrak{g})$ sending an automorphism to its derivative
at the unit element of $e$ is known to be an isomorphism of Lie groups.
Recall that the group $\mathrm{I}_G$ of all conjugations $g'\mapsto I_g(g'):=
gg'g^{-1}$ is a normal analytic subgroup of $\mathrm{Aut}_0(G)$ which is
isomorphic (by the above canonical map) to the adjoint Lie group, 
$\mathrm{Ad}_G$, of all adjoint representions
$\mathrm{Ad}_g$, $g\in G$, which is a normal analytic subgroup of 
$\mathrm{Aut}(\mathfrak{g})$.
For any strictly positive real number $\tau$, we set 
\begin{equation}\label{EqDefUTauInLieAlgebra}
 \mathcal{U}_{\tau\mathbf{i}}
   :=\{\xi\in \mathfrak{g}~|~\mathrm{all~the~eigenvalues~of~}\mathrm{ad}_\xi
               \mathrm{~lie~in~}S_{\tau\mathbf{i}}\},
\end{equation}
see eqn (\ref{EqDefOpenStripOfC}) for the definition of $S_{\tau\mathbf{i}}$,
and recall for any $\xi\in\mathfrak{g}$ the definition of its adjoint
representation $\mathrm{ad}_\xi:\eta\mapsto [\xi,\eta]$ for all 
$\eta\in\mathfrak{g}$. Furthermore, set 
\begin{equation}
  \mathcal{V}_{\tau\mathbf{i}}:= \exp\big(\mathcal{U}_{\tau\mathbf{i}}\big)  
           \subset G.
\end{equation}
We shall prove the following
\begin{prop}\label{PAdInvariantChartDomainExpAndBumpFunction}
 With the above definitions and notations, we have the following: Let
 $\tau$ be any real number such that $0<\tau\leq\pi$.
 \begin{enumerate}
  \item The subset 
    $\mathcal{U}_{\tau\mathbf{i}}$ is an open 
    $\mathrm{Aut}_0(\mathfrak{g})$-invariant
  neighbourhood of $0\in\mathfrak{g}$ such that for all 
  $\xi\in\mathcal{U}_{\tau\mathbf{i}}$ and for all $\eta$ in the nilradical 
  of $\mathfrak{g}$, the element $\xi+\eta$ still lies in 
    $\mathcal{U}_{\tau\mathbf{i}}$.
   \item The restriction of the exponential map to 
   $\mathcal{U}_{\tau\mathbf{i}}$ is a diffeomorphism onto
   $\mathcal{V}_{\tau\mathbf{i}}$ which is an open
  $\mathrm{Aut}_0(G)$-invariant neighbourhood of the unit element $e\in G$.
   \item Let $\tau'$ be any real number such that $0<\tau'<\tau\leq \pi$.
   Then there is a smooth $\mathrm{Aut}_0(\mathfrak{g})$-invariant 
   real-valued function
   $\gamma:\mathfrak{g}\to \mathbb{R}$ such that
   \begin{enumerate}
    \item $\gamma(\mathfrak{g})\subset [0,1]$,
    \item For all $\xi\in\mathcal{U}_{\tau'\mathbf{i}}$, we have
           $\gamma(\xi)=1$,
    \item The support of $\gamma$ is contained in 
         $\mathcal{U}_{\tau\mathbf{i}}$.
   \end{enumerate}
 \end{enumerate}
\end{prop}
\begin{prooof}
 \textbf{1.} Let $n:=\dim(\mathfrak{g})$ and define the map 
  $\tilde{\chi}:\mathfrak{g}\to \mathbb{C}[\lambda]^1_n$ for all 
  $\xi\in\mathfrak{g}$ by
  \[
      \tilde{\chi}(\xi)(\lambda):= \chi\big(\mathrm{ad}_{\xi}\big)(\lambda)
  \]
  where we have written $\chi$ for the characteristic polynomial.
  $\tilde{\chi}$ clearly is a polynomial map, hence a smooth map.
  Next, for all $\vartheta\in\mathrm{Aut}_0(\mathfrak{g})$ and all 
  $\xi\in\mathfrak{g}$, we get
   in $\mathcal{A}=\mathrm{Hom}_\mathbb{R}(\mathfrak{g},\mathfrak{g})$
   \begin{eqnarray*}
   \tilde{\chi}\big(\vartheta(\xi)\big)(\lambda)
   & = &
   \chi\big(\mathrm{ad}_{\vartheta(\xi)}\big)(\lambda) 
        =  \det\big(\lambda I - \mathrm{ad}_{\vartheta(\xi)}\big)
        = \det\big(\lambda I - 
           \vartheta\mathrm{ad}_{\xi}\vartheta^{-1}\big) \\
       & = & \det\big(\lambda I - \mathrm{ad}_{\xi}\big) = 
         \chi\big(\mathrm{ad}_{\xi}\big)(\lambda) 
         =\tilde{\chi}(\xi)(\lambda),
   \end{eqnarray*}
   whence $\tilde{\chi}$ is $\mathrm{Aut}_0(\mathfrak{g})$-invariant, 
   and thus the set
   \[
       \mathcal{U}_{\tau \mathbf{i}}=  
       \tilde{\chi}^{-1}\big(\mathcal{T}(S_{\tau\mathbf{i}}^{\times n})\big)
   \]
   is an open $\mathrm{Aut}_0(\mathfrak{g})$-invariant subset of
   $\mathfrak{g}$. Let 
   $\mathfrak{n}\subset \mathfrak{g}$ be the nilradical of $\mathfrak{g}$,
   i.e. the largest nilpotent ideal of $\mathfrak{g}$. Then it is well-known
   that there is a positive 
   integer $N$ such that for any $\eta_1,\ldots,\eta_N\in\mathfrak{n}$, the 
   product of linear maps $\mathrm{ad}_{\eta_1}\cdots\mathrm{ad}_{\eta_N}$
   vanishes. More generally, let $W_{N,m}$ let be any product of $N+m$
   adjoint representations $\mathrm{ad}_{\xi_i}$ with $1\leq i\leq N+m$
   and $\xi_i\in\mathfrak{g}$ for all $1\leq i\leq N+m$ such that $N$ adjoint
   representations are of the type $\mathrm{ad}_{\eta_j}$ with
   $\eta_1,\ldots,\eta_N\in \mathfrak{n}$. It is easy to see by induction
   on the positive integer $m$ that $W_{N,m}=0$: Indeed, this is clear for 
   $m=0$, and the induction step is simple if the first or the last position
   of $W_{N,m+1}$ is not in $\mathrm{ad}_\mathfrak{n}$, because then 
   $W_{N,m+1}$ contains
   $W_{N,m}$ as a factor and has to vanish. In the other cases, it is not 
   hard to 
   see that the $N$ $\mathrm{ad}_{\eta_j}$'s can all be moved to the right
   by getting commutators with the remaining $\mathrm{ad}_{\xi_i}$'s which 
   are adjoint representations of elements of the nilradical. Each commutator 
   reduces the word length by maintaining the number $N$ of elements in the 
   nilradical, and the term vanishes by induction. The remaining final term
   has a factor $\mathrm{ad}_{\eta_1}\cdots\mathrm{ad}_{\eta_N}$ and also has 
   to vanish, thus proving the induction. As a consequence, for any 
   $\xi_1,\ldots,\xi_m\in\mathfrak{g}$, any integer $1\leq s\leq m$ and
   for all $\eta\in\mathfrak{n}$,
   the linear map 
   $\mathrm{ad}_{\xi_1}\cdots \mathrm{ad}_{\xi_s}\mathrm{ad}_\eta
   \mathrm{ad}_{\xi_{s+1}}\cdots\mathrm{ad}_{\xi_m}$
   is nilpotent whence for all integers $1\leq r\leq n$
   \[
        \mathrm{trace}\Big(\big(\mathrm{ad}_{\xi+\eta}\big)^r\Big)=
         \mathrm{trace}\Big(\big(\mathrm{ad}_{\xi}\big)^r\Big)
   \]
   implying --thanks to the Waring identities, see e.g.~\cite[p.430]{AFI92}--
   \[
        \tilde{\chi}(\xi+\eta)(\lambda)=
          \det\big(\lambda I-\mathrm{ad}_\xi - \mathrm{ad}_\eta\big)
          =\det\big(\lambda I-\mathrm{ad}_\xi\big) =\tilde{\chi}(\xi)(\lambda)
   \]
   which shows the invariance of $\tilde{\chi}$ under translation by elements
   in the nilradical which proves the statement.\\
   \textbf{2.} For any $\xi,\eta\in\mathfrak{g}$ recall the formula for the 
   derivative of the exponential map:
   \[
       T_\xi\exp(\eta)=\left.\frac{d}{dt}\exp(\xi+t\eta)\right|_{t=0}
         = T_eL_{\exp(\xi)}
         \left(\frac{I-e^{\mathrm{ad}_\xi}}{\mathrm{ad}_\xi}(\eta)\right)  
         =T_eL_{\exp(\xi)}\big(h\big(\mathrm{ad}_\xi\big)(\eta)\big)
   \]
  where for any $g\in G$ $L_g:G\to G$ denotes the usual left multiplication 
  map in $G$, and we have used the function $h$,
  see eqn (\ref{EqDefFunctionHDerivativeOfExponential}).
  Since  $T_eL_{\exp(\xi)}:\mathfrak{g}\to T_{\exp(\xi)}G$ is
  a linear isomorphism, it remains to check the linear map 
  $h\big(\mathrm{ad}_\xi\big)$. But according to the second statement of
  Proposition 
  \ref{PInjectivityRegularityExpLinearMaps}, it follows that 
  $h\big(\mathrm{ad}_\xi\big)$ is a linear bijection for all $\xi\in
  \mathcal{U}_{\tau i}$. Hence the restriction of $\exp$ to
  $\mathcal{U}_{\tau i}$ is a local diffeomorphism. By the Inverse Function
  Theorem, the image $\mathcal{V}_{\tau i}=\exp\big(\mathcal{U}_{\tau i}\big)$
  must be an open set of $G$ containing $e$: Indeed, let 
  $g\in\mathcal{V}_{\tau i}$. Then there is $\xi\in \mathcal{U}_{\tau i}$
  with $g=\exp(\xi)$. There is an open neighbourhood of $\xi$ (which we can 
  choose to be in the open set $\mathcal{U}_{\tau i}$) on which the restriction
  of the exponential map is invertible, i.e. there is an open neighbourhood
  of $g$ and a local inverse of $\exp$. It follows that this second 
  neighbourhood is still contained in $\mathcal{V}_{\tau i}$ whence 
  $\mathcal{V}_{\tau i}$ is an open set of $G$ containing $e=\exp(0)$.
  Next, let $\xi,\eta\in \mathcal{U}_{\tau i}$ such that $\exp(\xi)=\exp(\eta)$:
  It follows that
  \[
      e^{\mathrm{ad}_\xi}=\mathrm{Ad}_{\exp(\xi)}
          =\mathrm{Ad}_{\exp(\eta)}=e^{\mathrm{ad}_\eta}.
  \]
  Again the third statement of Proposition 
  \ref{PInjectivityRegularityExpLinearMaps} implies that 
  \[
     \mathrm{ad}_\xi=\mathrm{ad_\eta}~~~\mathrm{hence}~~\exists~
      \zeta\in\mathfrak{z}(\mathfrak{g})~~\mathrm{such~that~}\eta=\xi+\zeta
  \]
  where $\mathfrak{z}(\mathfrak{g})=\mathrm{Ker}(\mathrm{ad})$ is the centre
  of $\mathfrak{g}$. It follows that
  \[
     \exp(\xi)=\exp(\eta)= \exp(\xi+\zeta)=\exp(\xi)\exp(\zeta)
  \]
  whence
  \[
        \exp(\zeta)=e.
  \]
  Since $G$ is a connected simply connected Lie group, $\zeta$ has to vanish by 
  the
  following beautiful argument of \cite[Section 2.3]{LT66}: Suppose that
  there is a nonvanishing $\zeta\in\mathfrak{z}(\mathfrak{g})$ such that
  $\exp(\zeta)=e$. It follows that 
  \[
       S:=\{\exp(t\zeta)\in G~|~t\in\mathbb{R}\}
  \]
  is a circle subgroup $S$ in the connected
  component $Z_0$ of the centre $Z$ of $G$. Consider the connected simply 
  connected Lie group $G'=G\times G$. There is a $2$-torus $S\times S$ in
  the connected component of the centre of $G'$. Let $\mathfrak{d}$ be a 
  one-dimensional subspace in the $2$-dimensional Lie algebra of $S\times S$
  such that its image under the exponential map is a \emph{dense} 
  one-dimensional
  central subgroup $D\subset S\times S$. Since $\mathfrak{d}$ is in the centre
  of $\mathfrak{g}'=\mathfrak{g}\oplus\mathfrak{g}$, the quotient
  Lie algebra $\mathfrak{g}'':=\mathfrak{g}'/\mathfrak{d}$ is well-defined. Let 
  $q:\mathfrak{g}'\to \mathfrak{g}''$ be the canonical projection which is a
  homomorphism of Lie algebras. Let $G''$ be a connected simply connected
  Lie group whose Lie algebra is $\mathfrak{g}''$. Then there is a unique 
  smooth 
  Lie group homomorphism $\kappa:G'\to G''$ such that its derivative at 
  the unit 
  element $(e,e)$ of $G'$ is equal to $q$. It follows that $\kappa$ is
  a surjective submersion. Now for all $\zeta'\in\mathfrak{d}$,
  we get 
  \[
       \kappa\big(\exp_{G'}(\zeta')\big)
          = \exp_{G''}\big(q(\zeta')\big)=\exp_{G''}\big(0\big)
         = e'',
  \]
  whence $D$ lies in the kernel of $\kappa$. The latter is a closed subgroup
  of $G'$, and contains thus the closure of $D$ which is equal to the torus
  $S\times S$. It follows that the kernel of $\kappa$ is at least 
  two-dimensional
  contradicting the fact that the dimension of $G''$ is equal to the dimension
  of $G'$ minus $1$. Therefore $\zeta=0$, and the restriction of $\exp$
  to $\mathcal{U}_{\tau\mathbf{i}}$ is injective, hence a diffeomorphism
  onto its image $\mathcal{V}_{\tau\mathbf{i}}$.
  Finally, let $\theta\in\mathrm{Aut}_0(G)$ and 
  $g\in \mathcal{V}_{\tau\mathbf{i}}$. Then there is
  a unique $\xi\in \mathcal{U}_{\tau\mathbf{i}}$ with $\exp(\xi)=g$. We get
  \[
     \theta(g)=\theta\big(\exp(\xi)\big)= \exp\big(T_e\theta(\xi)\big)\in 
       \mathcal{V}_{\tau\mathbf{i}}
  \]
  since $\mathcal{U}_{\tau\mathbf{i}}$ is invariant by the Lie algebra 
  automorphism $T_e\theta$
  whence $\mathcal{V}_{\tau\mathbf{i}}$ is invariant by Lie group automorphisms
  in $\mathrm{Aut}_0(G)$.\\
  \textbf{3.} Consider the following subsets of $\mathbb{C}[\lambda]^1_n$:
  \[
     \mathcal{W}_{\tau'\mathbf{i}}:=\mathcal{T}\big(S_{\tau'\pi}^{\times n}\big),~~
     \overline{\mathcal{W}}_{\tau'\mathbf{i}}
       :=\mathcal{T}\big(\overline{S_{\tau'\pi}}^{\times n}\big),~~    
        \mathcal{W}_{\tau\mathbf{i}}:=\mathcal{T}\big(S_{\tau\pi}^{\times n}\big)
  \]
  where $\overline{S_{\tau'\pi}}$ is the closure of $S_{\tau'\pi}$, i.e. 
  the closed strip of all those complex numbers
  whose imaginary part lies in the interval $[-\tau'\pi,\tau'\pi]$. The obvious 
  inclusions $S_{\tau'\pi}\subset \overline{S_{\tau'\pi}}\subset S_{\tau\pi}$
  imply the inclusions
  $\mathcal{W}_{\tau'\mathbf{i}}\subset 
   \overline{\mathcal{W}}_{\tau'\mathbf{i}}\subset \mathcal{W}_{\tau\mathbf{i}}$
   where $\mathcal{W}_{\tau'\mathbf{i}}$ and $\mathcal{W}_{\tau\mathbf{i}}$
   are open subsets and $\overline{\mathcal{W}}_{\tau'\mathbf{i}}$ is closed
   (in fact, the closure of $\mathcal{W}_{\tau'\mathbf{i}}$). It follows that
   the two open sets $\mathcal{W}_{\tau\mathbf{i}}$ and
   $\mathbb{C}[\lambda]^1_n\setminus \overline{\mathcal{W}}_{\tau'\mathbf{i}}$
   cover $\mathbb{C}[\lambda]^1_n$. Let $(\gamma', 1-\gamma')$ be a smooth 
   partition of unity subordinate to the open cover 
   $\big(\mathcal{W}_{\tau\mathbf{i}},
   \mathbb{C}[\lambda]^1_n\setminus 
   \overline{\mathcal{W}}_{\tau'\mathbf{i}}\big)$ of 
   $\mathbb{C}[\lambda]^1_n$.
   It follows that all the values of the smooth real-valued function $\gamma'$
   lie in the interval $[0,1]$, that its support is contained in 
   $\mathcal{W}_{\tau\mathbf{i}}$, and that the support of $1-\gamma'$ is 
   contained in 
   $\mathbb{C}[\lambda]^1_n\setminus 
   \overline{\mathcal{W}}_{\tau'\mathbf{i}}$, whence $\gamma'(f)=1$
   for all $f\in \overline{\mathcal{W}}_{\tau'\mathbf{i}}$, in particular
   for all $f\in \mathcal{W}_{\tau'\mathbf{i}}$. It follows that the 
   composed function $\gamma=\gamma'\circ \tilde{\chi}:\mathfrak{g}\to\mathbb{R}$
   has all the properties in the statement since clearly
   $\mathcal{U}_{\tau\pi}= \tilde{\chi}^{-1}\big(\mathcal{W}_{\tau\pi}\big)$,
   $\mathcal{U}_{\tau'\pi}= \tilde{\chi}^{-1}\big(\mathcal{W}_{\tau'\pi}\big)$,
   and $\tilde{\chi}$ is $\mathrm{Aut}_0(\mathfrak{g})$-invariant.
\end{prooof}

\appendix

\end{document}